\documentclass[12pt]{article}
\usepackage[cp1251]{inputenc}
\usepackage[russian]{babel}
\usepackage{amssymb}
\usepackage{amsmath}

\begin{document}

\begin{center}
\textbf{\large Estimates in Beurling--Helson type theorems.
Multidimensional case}\footnote{This is the author's preprint of
the paper published in \emph{Mathematical Notes}, \textbf{90}:3,
(2011), 373-384. The text below may slightly vary from its finally
published version.}
\end{center}

\begin{center}
Vladimir Lebedev
\end{center}

\begin{quotation}
{\small\textsc{Abstract.} We consider the spaces $A_p(\mathbb
T^m)$ of functions $f$ on the $m$ -dimensional torus $\mathbb T^m$
such that the sequence of the Fourier coefficients
$\widehat{f}=\{\widehat{f}(k), ~k \in \mathbb Z^m\}$ belongs to
$l^p(\mathbb Z^m), ~1\leq p<2$. The norm on $A_p(\mathbb T^m)$ is
defined by $\|f\|_{A_p(\mathbb T^m)}=\|\widehat{f}\|_{l^p(\mathbb
Z^m)}$. We study the rate of growth of the norms
$\|e^{i\lambda\varphi}\|_{A_p(\mathbb T^m)}$ as
$|\lambda|\rightarrow \infty, ~\lambda\in\mathbb R,$ for $C^1$
-smooth real functions $\varphi$ on $\mathbb T^m$ (the
one-dimensional case was investigated by the author earlier). The
lower estimates that we obtain have direct analogues for the
spaces $A_p(\mathbb R^m)$.

  References: 15 items.


  AMS 2010 Mathematics Subject Classification 42B05, 42B35}
\end{quotation}

\quad

\begin{center}
\textbf{Introduction}
\end{center}

  Given any integrable function $f$ on the $m$ -dimensional torus
$\mathbb T^m=\mathbb R^m/2\pi\mathbb Z^m, ~m\geq 1$ (where
$\mathbb R$ is the real line, $\mathbb Z$ is the set of integers)
consider its Fourier coefficients:
$$
\widehat{f}(k)=\frac{1}{(2\pi)^m}\int_{\mathbb T^m} f(t)e^{-i(k, t)} dt,
\qquad k\in\mathbb Z^m.
$$

  Let $A_1(\mathbb T^m)$ be the space of continuous
functions $f$ on $\mathbb T^m$ such that the sequence of Fourier
coefficients $\widehat{f}=\{\widehat{f}(k), ~k\in\mathbb Z^m\}$
belongs to $l^1(\mathbb Z^m)$. For $1<p\leq 2$ let $A_p(\mathbb
T^m)$ be the space of integrable functions $f$ on $\mathbb T^m$
such that $\widehat{f}\in l^p(\mathbb Z^m)$. Provided with the
natural norms
$$
\|f\|_{A_p(\mathbb T^m)}=\|\widehat{f}\|_{l^p(\mathbb Z^m)}=
\bigg(\sum_{k\in\mathbb Z^m}|\widehat{f}(k)|^p\bigg)^{1/p}
$$
the spaces $A_p$ are Banach spaces $(1\leq p\leq 2)$. The space
$A=A_1$ is a Banach algebra (with the usual multiplication of
functions).

  According to the Beurling--Helson theorem [1] (see also [2]),
if $\varphi$ is a map of the circle $\mathbb T$ into itself such
that $\|e^{in\varphi}\|_{A(\mathbb T)}=O(1), ~n\in\mathbb Z$, then
$\varphi$ is linear (with integer tangent coefficient), i.e.
$\varphi(t)=kt+\varphi(0), ~k\in\mathbb Z$. A similar statement
holds for the maps $\varphi:\mathbb T^m\rightarrow\mathbb T$. This
case easily reduces to the one-dimensional case.

  Let $C^\nu (\mathbb T^m)$ be the class of (complex-valued)
functions on the torus $\mathbb T^m$ such that all partial
derivatives of order $\nu$ are continuous.

  In the present paper we study the growth of the norms
$\|e^{i\lambda\varphi}\|_{A_p(\mathbb T^m)}$ as
$\lambda\rightarrow\infty, ~\lambda\in\mathbb R,$ for $C^1$
-smooth real functions $\varphi$ on $\mathbb T^m$. In the
one-dimensional case we studied this question in [3]. The same
paper contains a survey on the subject.

  It is easy to show that for every $C^1$ -smooth (real) function
 $\varphi$ on the circle $\mathbb T$ we have
\footnote{Actually estimate (1) holds even in the case when
$\varphi$ is absolutely continuous with derivative in $L^2(\mathbb
T)$ and, in particular, when $\varphi$ satisfies the Lipschitz
condition of order $1$ (see [2, Ch. VI, \S~3] for $p=1$; for
$1<p<2$ the estimate follows immediately by interpolation between
$l^1$ and $l^2$).}
$$
\|e^{i\lambda\varphi}\|_{A_p(\mathbb T)}=
O(|\lambda|^{\frac{1}{p}-\frac{1}{2}}),
\qquad |\lambda|\rightarrow\infty, \quad \lambda\in\mathbb R,
\eqno(1)
$$
for all $p, ~1\leq p<2$.

   On the other hand we have the Leibenson--Kahane--Alp\'ar
estimate ([4], [5], [2], [6]): if $\varphi\in C^2(\mathbb T)$ is a
nonconstant (which is equivalent, due to periodicity, to it being
nonlinear) real function, then
$$
\|e^{i\lambda\varphi}\|_{A_p(\mathbb T)}\geq
c_p |\lambda|^{\frac{1}{p}-\frac{1}{2}}, \qquad \lambda\in\mathbb R,
\eqno(2)
$$
for all $1\leq p<2$.

   Thus, for every $C^2$ -smooth nonlinear real function
$\varphi$ on $\mathbb T$ we have \footnote{We write
$a(\lambda)\simeq b(\lambda)$ in the case when $c_1\leq
a(\lambda)/b(\lambda)\leq c_2$ for all sufficiently large
$|\lambda|$ (with constants $c_1, c_2>0$ independent of
$\lambda$).}
$$
\|e^{i\lambda\varphi}\|_{A_p(\mathbb T)}\simeq
|\lambda|^{\frac{1}{p}-\frac{1}{2}}.
\eqno(3)
$$
In particular $\|e^{i\lambda\varphi}\|_{A(\mathbb T)}\simeq
|\lambda|^{1/2}$.

  The case of $C^1$ -smooth phase $\varphi$ is essentially
different from the $C^2$ -smooth case. As we showed in [3] (and
earlier in [7]), for $\varphi\in C^1(\mathbb T)$ the norms
$\|e^{i\lambda\varphi}\|_{A(\mathbb T)}$ can grow rather slowly,
namely, one can closely approach the $O(\log |\lambda|)$
-condition. If $p>1$ the corresponding norms can even be bounded.

   We note that the proof of the Leibenson--Kahane--Alp\'ar
estimate (2) is based on the van der Corput lemma and essentially
uses nondegeneration of the curvature of a certain arc of the
graph of $\varphi$. This approach does not allow to consider the
functions of smoothness less then $C^2$.

  In the multidimensional case for the phase functions
$\varphi$ of smoothness $C^2$ (and higher) the behavior of the
norms $\|e^{i\lambda\varphi}\|_{A}$ was considered by Hedstrom
[8]. As in the one-dimensional case it is easy to get an upper
estimate; for instance, \footnote{To ensure this it suffices to
repeat with obvious modifications the arguments used for $m=1$ in
[2, Ch. VI, \S~3].} if $\varphi\in C^\nu (\mathbb T^m), ~\nu>m/2,
~m\geq 2,$ then $\|e^{i\lambda\varphi}\|_{A(\mathbb
T^m)}=O(|\lambda|^{m/2})$ (in [8] this estimate is obtained under
somewhat different assumptions on smoothness). The same work [8]
contains the following lower estimate: if $\varphi\in C^2(\mathbb
T^m)$ is a real function such that the determinant of the matrix
of its second derivatives is not identically equal to zero, then
$$
\|e^{i\lambda\varphi}\|_{A(\mathbb T^m)}\geq c |\lambda|^{m/2}.
\eqno(4)
$$
This is proved by reduction to the one-dimensional case.

  In \S~1 we obtain Theorem 1, which is the main result of
the present paper. In this theorem we give lower estimates for the
norms $\|e^{i\lambda\varphi}\|_{A_p(\mathbb T^m)}$ for $C^1$
-smooth real functions $\varphi$ on the torus $\mathbb T^m$. We
assume that the range $\nabla\varphi(\mathbb T^m)$ of the gradient
$\nabla\varphi$ of $\varphi$ has positive (Lebesgue) measure; in
this case we say that the gradient of $\varphi$ is nondegenerate.
The proof of the theorem is based on the natural modification of
the method that we used earlier for the one-dimensional case in
[3]. It can be called the concentration of Fourier transform large
values method. We also obtain local versions of the lower
estimates in the spaces $A_p(\mathbb T^m)$ (Theorem 1$'$) and
$A_p(\mathbb R^m)$ (Theorem 1$''$), where $A_p(\mathbb R^m)$ is
the space of functions $f$ on $\mathbb R^m$ such that the Fourier
transform $\widehat{f}$ belongs to $L^p(\mathbb R^m)$.

  We note that in the one-dimensional case we have
the obvious inclusion $C^1(\mathbb T)\subseteq A(\mathbb
T)\subseteq A_p(\mathbb T)$. At the same time for the torus of
dimension $m\geq 3$ the condition of $C^1$ -smoothness does not
guarantee that a function belongs to all classes $A_p, 1\leq p<2$.
(Smoothness conditions that depend on dimension and imply that a
function belongs to the classes $A_p(\mathbb T^m)$ are well known,
see \S~3). This is why we naturaly consider the lower estimates of
the norms $\|e^{i\lambda\varphi}\|_{A_p(\mathbb T^m)}$ for
$\lambda\in\Lambda(\varphi, p)$, where for a function $\varphi$ on
$\mathbb T^m$ the set $\Lambda(\varphi, p)$ is the set of all
those $\lambda\in\mathbb R$ for which $e^{i\lambda\varphi}\in
A_p(\mathbb T^m)$.

  Here, in Introduction, we note only the corollaries
of Theorem 1. Let $\varphi\in C^1(\mathbb T^m)$ be a real function
such that its gradient is nondegenerate and satisfies the
Lipschitz condition of order  $\alpha, ~0<\alpha\leq 1$; then
(Corollary 1) for $1\leq p<1+\alpha$ we have
$$
\|e^{i\lambda\varphi}\|_{A_p(\mathbb T^m)}\geq
c_p |\lambda|^{m\big(\frac{1}{p}-\frac{1}{1+\alpha}\big)},
\qquad \lambda\in \Lambda(\varphi, p).
$$
In particular (Corollary 2), we see that if the gradient of a
phase function $\varphi$ is nondegenerate and satisfies the
Lipschitz condition of order $1$, then
$$
\|e^{i\lambda\varphi}\|_{A_p(\mathbb T^m)}\geq c_p
|\lambda|^{m\big(\frac{1}{p}-\frac{1}{2}\big)},
\qquad \lambda\in \Lambda(\varphi, p).
\eqno(5)
$$
Putting here $p=1$ we see that every such function $\varphi$
satisfies (4) for all $\lambda\in\Lambda(\varphi, 1)$.

  In \S~2 for each class $C^{1, \omega}(\mathbb T^m)$
(see the definition of the classes $C^{1, \omega}$ at the end of
the introduction) we construct a real function $\varphi\in C^{1,
\omega}(\mathbb T^m)$ such that it has nondegenerate gradient and
the norms $\|e^{i\lambda\varphi}\|_{A_p(\mathbb T^m)}$ grow very
slowly. (For the one-dimensional case we did this in [3]. The
general case can be easily deduced from the one-dimensional case.)
Thus we show (see Theorem 2 and Corollaries 3 and 4) that the
lower estimates obtained in \S~1 are close to being sharp and in
certain cases are sharp.

  In \S~3 using quite standard methods
we obtain the multidimensional version of estimate (1), namely
(Theorem 3): if $\varphi$ is a sufficiently smooth (depending on
dimension) real function on $\mathbb T^m$, then
$\|e^{i\lambda\varphi}\|_{A_p(\mathbb
T^m)}=O(|\lambda|^{m(1/p-1/2)})$. Hence taking our estimate (5)
into account, we obtain the multidimensional version of relation
(3), namely (Theorem 4): if $\varphi$ is sufficiently smooth and
its gradient is nondegenerate, then
$$
\|e^{i\lambda\varphi}\|_{A_p(\mathbb T^m)}\simeq
|\lambda|^{m(1/p-1/2)},
$$
in particular $\|e^{i\lambda\varphi}\|_{A(\mathbb T^m)}\simeq
|\lambda|^{m/2}$.

\quad

  We use the following notation. Let $V$ be a domain in
$\mathbb R^m$ and let $g$ be a function on $V$. We define the
modulus of continuity of $g$ by
$$
\omega (V, g, \delta)=\sup_{\underset{t_1, t_2\in V}{|t_1-t_2|\leq
\delta}}|g(t_1)-g(t_2)|, \qquad \delta\geq 0,
$$
where $|x|$ is the length of a vector $x\in\mathbb R^m$. If
$V=\mathbb R^m$, then we just write $\omega (g, \delta)$. Let
$\omega$ be a given continuous nondecreasing function on $[0,
+\infty), ~\omega(0)=0$. The class $\mathrm{Lip_\omega}(V)$
consists of functions $g$ on $V$ satisfying $\omega (V, g,
\delta)=O(\omega(\delta)), ~\delta\rightarrow +0$. The class
$C^{1, \omega}(V)$ consists of functions $f$ on $V$ such that all
derivatives of the first order $\partial f/\partial t_j, ~j=1, 2,
\ldots, m,$ belong to $\mathrm{Lip_\omega}(V)$. Certainly for a
(real) function $\varphi$ on $V$ the condition $\varphi\in C^{1,
\omega}(V)$ means that $\omega (V, \nabla\varphi,
\delta)=O(\omega(\delta)), ~\delta\rightarrow +0$, where
$$
\omega (V, \nabla \varphi, \delta)=
\sup_{\underset{t_1, t_2\in V}{|t_1-t_2|\leq
\delta}}|\nabla \varphi(t_1)-\nabla \varphi(t_2)|, \qquad \delta\geq 0,
$$
is the modulus of continuity of the gradient $\nabla\varphi$ of
$\varphi$. The class $C^{1, \omega}(\mathbb T^m)$ consists of
functions that are  $2\pi$ -periodic with respect to each variable
and belong to $C^{1, \omega}(\mathbb R^m)$. For $0<\alpha\leq 1$
we write $C^{1, \alpha}$ instead of $C^{1, \delta^\alpha}$.
Generally, for an arbitrary $\nu=0, 1, 2, \ldots$ and
$0<\alpha\leq 1$ let $C^{\nu, \alpha}(\mathbb T^m)$ be the class
of functions $f$ on $\mathbb T^m$ such that $f$ is $\nu$ times
differentiable and all partial derivatives of order $\nu$ of $f$
(for $\nu=0$ the function $f$ itself) satisfy the Lipschitz
condition of order $\alpha$ (i.e. belong to
$\mathrm{Lip_{\delta^\alpha}}$). It is convenient to put $C^{\nu,
0}=C^\nu$. For an arbitrary measurable set $E$ in $\mathbb T^m$ or
in $\mathbb R^m$ by $|E|$ we denote its Lebesgue measure. By $(x,
y)$ we denote the usual inner product of vectors $x$ and $y$ in
$\mathbb R^m$ (or of vectors $x\in \mathbb Z^m, ~y\in \mathbb
T^m$). If $W$ is a set in $\mathbb R^m$ and $\lambda\in \mathbb
R$, then we put $\lambda W=\{\lambda x: x\in W\}$. In the usual
way we identify integrable functions on the torus $\mathbb T^m$
with integrable functions on the cube $[0, 2\pi]^m$. By $c, c_p,
c(p, \varphi), c_m,$ etc. we denote various positive constants
which may depend only on  $p, \varphi$ and the dimension $m$.

\quad

\begin{center}
\textbf{\S~1. Lower estimates}
\end{center}

   As we indicated in Introduction, given a function $\varphi$
on the torus $\mathbb T^m$ we denote by $\Lambda(\varphi, p)$ the
set of all those $\lambda\in\mathbb R$ for which
$e^{i\lambda\varphi}\in A_p(\mathbb T^m)$.

\quad

\textbf{Theorem 1.} \emph{Let $1\leq p <2$. Let $\varphi\in C^{1,
\omega}(\mathbb T^m)$ be a real function. Suppose that the
gradient $\nabla\varphi$ is nondegenerate that is the set
$\nabla\varphi (\mathbb T^m)$ is of positive measure. Then
$$
\|e^{i\lambda\varphi}\|_{A_p(\mathbb T^m)} \geq
c \bigg(|\lambda|^{1/p}\chi^{-1}\bigg(\frac{1}{|\lambda|}\bigg)\bigg)^m,
\qquad \lambda\in\Lambda(\varphi, p), \quad |\lambda|\geq 1,
$$
where $\chi^{-1}$ is the function inverse to
$\chi(\delta)=\delta\omega(\delta)$ and $c=c(p, \varphi)>0$ is
independent of $\lambda$.}

\quad

  In the one-dimensional case we obtained this theorem in
[3]. (For each $C^1$ -smooth function $\varphi$ on $\mathbb T$ we
have $\Lambda(\varphi, p)=\mathbb R$ for all $p\geq 1$.
Nondegeneration of the gradient in the one-dimensional case means
nonlinearity of $\varphi$ which due to periodicity is equivalent
to the condition that $\varphi\neq\mathrm{const}$.)

  Theorem 1 immediately implies the following corollary.

\quad

\textbf{Corollary 1.} \emph{Let $0<\alpha\leq 1$. Let $\varphi\in
C^{1, \alpha}(\mathbb T^m)$ be a real function with nondegenerate
gradient. Then for all $p, ~1\leq p<1+\alpha,$ we have
$$
\|e^{i\lambda\varphi}\|_{A_p(\mathbb T^m)}\geq
c_p |\lambda|^{m\big(\frac{1}{p}-\frac{1}{1+\alpha}\big)},
\qquad \lambda\in\Lambda(\varphi, p).
$$
In particular $\|e^{i\lambda\varphi}\|_{A(\mathbb T^m)}\geq c
|\lambda|^{\frac{m\alpha}{1+\alpha}}, ~\lambda\in \Lambda(\varphi,
1)$.}

\quad

  We especially note the case of a $C^2$ -smooth phase and even the
more general case of a $C^{1, 1}$ -smooth phase.

\quad

\textbf{Corollary 2.} \emph{Let $\varphi\in C^{1, 1}(\mathbb T^m)$
be a real function with nondegenerate gradient. Then for all $p,
~1\leq p<2,$ we have
$$
\|e^{i\lambda\varphi}\|_{A_p(\mathbb T^m)}\geq
c_p |\lambda|^{m\big(\frac{1}{p}-\frac{1}{2}\big)},
\qquad \lambda\in\Lambda(\varphi, p).
$$
In particular $\|e^{i\lambda\varphi}\|_{A(\mathbb T^m)}\geq c
|\lambda|^{m/2}, ~\lambda\in\Lambda(\varphi, 1)$.}

\quad

  We shall see that a local version of Theorem 1 also holds. Let
$E$ be an arbitrary set contained in $[0, 2\pi]^m$ or, more
generally, contained in some cube with edges of length  $2\pi$
parallel to coordinate axes. We say that a function $f$ defined on
$E$ belongs to $A_p(\mathbb T^m, E)$ if there exists a function
$F\in A_p(\mathbb T^m)$ such that its restriction $F_{|E}$ to the
set $E$ coincides with $f$. We put
$$
\|f\|_{A_p(\mathbb T^m, E)}=\inf_{F_{|E}=f}\|F\|_{A_p(\mathbb T^m)}.
$$

   As in Theorem 1, everywhere below $\chi^{-1}$ is the function
inverse to $\chi(\delta)=\delta\omega(\delta)$.

\quad

    \textbf{Theorem $1'$.} \emph{Let $1\leq p <2$. Let $V$ be a
domain in $[0, 2\pi]^m$. Let $\varphi\in C^{1, \,\omega}(V)$ be a
real function. Suppose that the gradient $\nabla\varphi$ is
nondegenerate on $V$, that is, the set $\nabla\varphi (V)$ is of
positive measure. Then
$$
\|e^{i\lambda\varphi}\|_{A_p(\mathbb T^m, V)} \geq
c \bigg(|\lambda|^{1/p}\chi^{-1}\bigg(\frac{1}{|\lambda|}\bigg)\bigg)^m
$$
for all those $\lambda\in\mathbb R, ~|\lambda|\geq 1,$ for which
$e^{i\lambda\varphi}\in A_p(\mathbb T^m, V)$.}

\quad

  The local version of Corollary 1 (as well as of Corollary 2)
is obvious.

  Technically it will be convenient to deal with the
spaces $A_p$ in nonperiodic case. (The reason for this is that in
distinction with the one-dimensional case, for a $C^1$ -smooth
function of several variables the range of its gradient can be
very complicated. \footnote{In this connection we note the
following question. Let $V$ be a domain in $\mathbb R^m$ and let
$\varphi\in C^1(V)$ be a real function. Suppose that the set
$\nabla \varphi (V)$ is of positive measure. Is it true then that
this set has nonempty interior? For $m=2$ the answer to this
question is positive [9]. For $m\geq 3$ the answer is unknown.})

   Let $A_p(\mathbb R^m)$, where $1\leq p\leq\infty$, be the space
of tempered distributions $f$ on $\mathbb R^m$ such that the
Fourier transform $\widehat{f}$ belongs to $L^p(\mathbb R^m)$. We
put
$$
\|f\|_{A_p(\mathbb R^m)}=\|\widehat{f}\|_{L^p(\mathbb R^m)}=
\bigg(\int_{\mathbb R^m} |\widehat{f}(u)|^p du \bigg)^{1/p}.
$$
For $1\leq p\leq 2$ each distribution that belongs to $A_p$ is
actually a function in $L^q, ~1/p+1/q=1$. For $p=1$ we naturally
assume that $f$ is continuous.\footnote{We note that in [1] the
notation $A$ stands for the space of (inverse) Fourier transforms
of measures on $\mathbb R$. We follow the notation which is common
nowadays (see e.g. [2]).}

  We choose the normalization factor of the Fourier transform so
that
$$
\widehat{f}(u)=\frac{1}{(2\pi)^m}\int_{\mathbb R^m} f(t)e^{-i(u, t)} dt,
\qquad u\in\mathbb R^m,
$$
for $f\in L^1(\mathbb R^m)$.

  We often write $A$ instead of $A_1$.

  Let us define the local spaces $A_p$ in nonperiodic case.

  Let $E\subseteq \mathbb R^m$ be an arbitrary set. We say that a function
$f$ defined on $E$ belongs to $A_p (\mathbb R^m, E)$ if there
exists a function $F\in A_p(\mathbb R^m)$ such that its
restriction $F_{|E}$ to the set $E$ coincides with $f$. The norm
on $A_p (\mathbb R^m, E)$ is defined in the natural way:
$$
\|f\|_{A_p(\mathbb R^m, E)}=\inf_{F_{|E}=f} \|F\|_{A_p(\mathbb R^m)}.
$$

\quad

    \textbf{Theorem $1''$.} \emph{Let $1\leq p <2$. Let $V$ be a
domain in $\mathbb R^m$. Let $\varphi\in C^{1, \,\omega}(V)$ be a
real function on $V$. Suppose that the gradient $\nabla\varphi$ is
nondegenerate on $V$ that is the set $\nabla\varphi (V)$ is of
positive measure. Then
$$
\|e^{i\lambda\varphi}\|_{A_p(\mathbb R^m, V)} \geq
c \bigg(|\lambda|^{1/p}\chi^{-1}\bigg(\frac{1}{|\lambda|}\bigg)\bigg)^m
$$
for all those $\lambda\in\mathbb R, ~|\lambda|\geq 1,$ for which
$e^{i\lambda\varphi}\in A_p(\mathbb R^m, V)$.}

\quad

  Theorem 1, which is the main result of this section,
follows immediately from Theorem  $1'$. Let us show that in turn
Theorem $1'$ follows from Theorem $1''$; then we shall prove
Theorem $1''$ itself.

  Fix $p, ~1\leq p<2$. For each $\lambda$ such that
$e^{i\lambda\varphi}\in A_p(\mathbb T^m, V)$ consider an arbitrary
$2\pi$ -periodic (with respect to each variable) extension
$F_\lambda\in A_p(\mathbb T^m)$ of the function
$e^{i\lambda\varphi}$ from $V$ to $\mathbb R^m$.

  Under the assumptions of Theorem $1'$ we can find a cube
$I\subseteq V$ with edges parallel to the coordinate axes such
that together with its closure it is contained in the interior of
the cube $[0, 2\pi]^m$ and the gradient $\nabla\varphi$ is
nondegenerate on $I$. Let $\chi$ be an infinitely differentiable
function on $\mathbb R^m$ equal to $1$ on $I$ and equal to $0$ on
the compliment $\mathbb R^m\setminus[0, 2\pi]^m$. We have $\chi\in
A(\mathbb R^m)$. Put $c_0=\|\chi\|_{A(\mathbb R^m)}$. For each
$u\in\mathbb R^m$ define the function $e_u$ on $\mathbb R^m$ by
$e_u(t)=e^{i(u, t)}, ~t\in\mathbb R^m$. We have
$\|e_u\chi\|_{A(\mathbb R^m)}=c_0$.

  For an arbitrary function $h$ on $\mathbb R^m$
that vanishes outside of the cube $[0, 2\pi]^m$ let
$\widetilde{h}$ be its $2\pi$ -periodic with respect to each
variable extension from $[0, 2\pi]^m$ to $\mathbb R^m$. It is well
known (see e.g. [10, \S~50]) that if $h\in A(\mathbb R^m)$ then
$\widetilde{h}\in A(\mathbb T^m)$ and
$\|\widetilde{h}\|_{A(\mathbb T^m)}\leq c\|h\|_{A(\mathbb R^m)}$.

  We put $g_u=\widetilde{ e_u\chi}$. For each $u\in\mathbb R^m$
we have $\|g_u\|_{A(\mathbb T^m)}\leq c_1$ whence $\|g_u
F_\lambda\|_{A_p(\mathbb T^m)}\leq c_1 \|F_\lambda\|_{A_p(\mathbb
T^m)}$. It is also clear that $\|g_u F_\lambda\|_{A_p(\mathbb
T^m)}$ is a measurable function of $u$.

  The function $\chi F_\lambda$ coincides with $e^{i\lambda\varphi}$
on $I$. Thus,
$$
\|e^{i\lambda\varphi}\|_{A_p(\mathbb R^m, I)}^p
\leq\|\chi F_\lambda\|_{A_p(\mathbb R)}^p
=\sum_{k\in\mathbb Z^m} \int_{[0, 1]^m}
|\widehat{\chi F_\lambda}(u+k)|^p du
$$
$$
=\int_{[0, 1]^m}\sum_{k\in\mathbb Z^m} |\widehat{e_{-u}\chi
F_\lambda}(k)|^p du=
\int_{[0, 1]^m} \|g_{-u}F_\lambda\|_{A_p(\mathbb T^m)}^p du\leq
c_1 \|F_\lambda\|_{A_p(\mathbb T^m)}^p.
$$
It remains to use Theorem 1$''$.

\quad

  \emph{Proof of Theorem} $1''$. We can find a closed cube
$I\subseteq V$ with edges parallel to the coordinate axes such
that its image $W=\nabla\varphi(I)$ is of positive measure. The
set $W$ is bounded.

  Fix $c>0$ so that
$$
\omega(I,\nabla \varphi, \delta)\leq c \omega(\delta), \qquad \delta\geq
0.
$$

    For each $\lambda>0$ choose $\delta_\lambda>0$ so that
$$
\chi(m^{1/2}2\delta_\lambda)=\frac{1}{2c\lambda}.
\eqno(6)
$$

  For $\varepsilon>0$ let $\Delta _\varepsilon$
mean the ``triangle'' function supported on the interval
$(-\varepsilon, \varepsilon)$, that is the function on $\mathbb R$
defined by
$$
\Delta _\varepsilon (t)=\max \Big (1-\frac{|t|}{\varepsilon },~0 \Big ),
\quad t \in \mathbb R,
$$
and for an arbitrary interval $J\subseteq\mathbb R$ let $\Delta
_J$ be the triangle function supported on $J$ that is $\Delta
_J(t)=\Delta _{|J|/2}(t-c_J)$, where $c_J$ is the center of the
interval $J$ (and $|J|$ is its length). Let then $J$ be a cube in
$\mathbb R^m$ with edges parallel to coordinate axes, $J=J_1\times
J_2 \times \ldots \times J_m$. We define the triangle function
$\Delta_J$ supported on $J$ as follows
$$
\Delta_J(t)=\Delta_{J_1}(t_1)\Delta_{J_2}(t_2)\ldots \Delta_{J_m}(t_m),
\qquad t=(t_1, t_2, \ldots, t_m)\in\mathbb R^m.
$$

  We shall use the following lemma.

\quad

\textbf{Lemma 1.} \emph{Let $\lambda>0$ be sufficiently large. Let
$F_\lambda$ be an arbitrary function on $\mathbb R^m$ that
coincides with $e^{i\lambda\varphi}$ on $V$. Then for each
$u\in\lambda W$ there exists a cube $I_{\lambda, u}\subseteq I$
with edges of length $2\delta_\lambda$ parallel to coordinate
axes, such that
$$
|(\Delta_{I_{\lambda, u}}F_\lambda)^{^\wedge}(u)|\geq
c_m\delta_\lambda^m,
$$
where $c_m>0$ depends only on the dimension $m$.}

\quad

\emph{Proof}. Denote by $a$ the length of the edge of the cube
$I$. We shall assume that $\lambda>0$ is so large that
$$
2\delta_\lambda<a.
\eqno(7)
$$

   Take an arbitrary $u\in\lambda W$. We can find a point
$t_{\lambda, u}\in I$ such that $\nabla\varphi(t_{\lambda,
u})=\lambda^{-1}u$. Let $I_{\lambda, u}\subseteq I$ be a cube (a
closed one if necessary) with edges of length $2\delta_\lambda$
parallel to coordinate axes such that it contains the point
$t_{\lambda, u}$ (see (7)). Consider the following linear
function:
$$
\varphi_{\lambda, u}(t)=\varphi(t_{\lambda, u})+(\lambda^{-1}u,
~t-t_{\lambda, u}), \qquad t\in\mathbb R^m.
$$

   If $t\in I_{\lambda, u}$ then for some point
$\theta\in I_{\lambda, u}$ we shall have
$$
\varphi(t)-\varphi(t_{\lambda, u})=(\nabla \varphi(\theta),
~t-t_{\lambda, u})
$$
($\theta$ lies on the strait segment that joins  $t$ with
$t_{\lambda, u}$) and therefore
$$
|\varphi(t)-\varphi_{\lambda, u}(t)|=|\varphi(t)-\varphi(t_{\lambda,
u})-(\nabla \varphi(t_{\lambda, u}), ~t-t_{\lambda, u})|
$$
$$
=|(\nabla \varphi(\theta)-\nabla \varphi(t_{\lambda, u}), ~t-t_{\lambda,
u})|\leq |\nabla \varphi(\theta)-\nabla \varphi(t_{\lambda, u})|
|t-t_{\lambda, u}|
$$
$$
\leq \omega(I,\nabla\varphi,
m^{1/2}2\delta_\lambda)m^{1/2}2\delta_\lambda\leq c\chi (m^{1/2}
2\delta_\lambda).
$$
Hence, taking into account (6), we see that
$$
|e^{i\lambda\varphi(t)}-e^{i\lambda\varphi_{\lambda, u}(t)}|\leq
|\lambda\varphi(t)-\lambda\varphi_{\lambda, u}(t)|\leq \lambda c\chi
(m^{1/2} 2\delta_\lambda)=\frac{1}{2}, \qquad t\in I_{\lambda, u}.
$$

   Using this estimate we obtain
$$
|(\Delta_{I_{\lambda, u}} F_\lambda)^{^\wedge}(u)-
(\Delta_{I_{\lambda,
u}}e^{i\lambda\varphi_{\lambda, u}})^{^\wedge}(u)|
\leq\frac{1}{(2\pi)^m} \int_{I_{\lambda, u}} \Delta_{I_{\lambda,
u}}(t)|e^{i\lambda\varphi(t)}-e^{i\lambda\varphi_{\lambda, u}(t)}|
dt
$$
$$
\leq \frac{1}{2}\cdot\frac{1}{(2\pi)^m}
\int_{\mathbb R^m}\Delta_{I_{\lambda,
u}}(t)dt=\frac{1}{2}\widehat{\Delta_{I_{\lambda, u}}}(0).
$$
At the same time
$$
|(\Delta_{I_{\lambda, u}}e^{i\lambda\varphi_{\lambda,
u}})^{^\wedge}(u)|=\bigg |\frac{1}{(2\pi)^m} \int_{\mathbb R^m}
\Delta_{I_{\lambda, u}}(t)e^{i(\lambda\varphi_{\lambda, u}(t)-(u, t))}
dt\bigg |=\widehat{\Delta_{I_{\lambda, u}}}(0).
$$
Thus
$$
|(\Delta_{I_{\lambda, u}}F_\lambda)^{^\wedge}(u)|\geq
\frac{1}{2}\widehat{\Delta_{I_{\lambda, u}}}(0)=c_m \delta_\lambda^m.
$$
The lemma is proved.

\quad

    Fix $p, ~1\leq p<2$. Everywhere below we assume that
the frequencies $\lambda$ are such that $e^{i\lambda\varphi}\in
A_p(\mathbb R^m, V)$. We can assume that the set of these
$\lambda$ is unbounded, otherwise there is nothing to prove. We
can also assume that these frequencies $\lambda$ are positive (the
complex conjugation does not affect the norm of a function in
$A_p$). For every such $\lambda$ let $F_\lambda$ be an extension
of the function $e^{i\lambda\varphi}$ from $V$ to $\mathbb R^m$
satisfying
$$
\|F_\lambda\|_{A_p(\mathbb R^m)}\leq
2\|e^{i\lambda\varphi}\|_{A_p(\mathbb R^m, V)}.
\eqno(8)
$$
Let $I_{\lambda, u}$ be the corresponding cubes whose existence
(for all sufficiently large $\lambda$) is established in Lemma 1.

  Define functions $g_\lambda$ by
$$
g_\lambda=(|(F_\lambda)^{^\wedge}|)^{^\vee},
$$
where $^\vee$ means the inverse Fourier transform. We have
$g_\lambda\in A_p(\mathbb R^m)$.

   It is well known that the function $\Delta_\varepsilon$ belongs
to $A(\mathbb R)$ and has nonnegative Fourier transform. Hence
$\widehat{\Delta_{(-\varepsilon, \varepsilon)^m}}\geq 0$ and since
for an arbitrary cube $J\subseteq \mathbb R^m$ with edges of
length $2\varepsilon$ parallel to coordinate axes the function
$\Delta_J$ is obtained from $\Delta_{(-\varepsilon,
\varepsilon)^m}$ by shift, we have
$|\widehat{\Delta_J}|=\widehat{\Delta_{(-\varepsilon,
\varepsilon)^m}}$. Therefore ($\ast$ means convolution),
$$
|(\Delta_{I_{\lambda,
u}}F_\lambda)^{^\wedge}(u)|
=|(\Delta_{I_{\lambda, u}})^{^\wedge}\ast
(F_\lambda)^{^\wedge} (u)|\leq |(\Delta_{I_{\lambda,
u}})^{^\wedge}|\ast |(F_\lambda)^{^\wedge}|(u)=
$$
$$
= (\Delta_{(-\delta_\lambda, ~\delta_\lambda)^m})^{^\wedge}\ast
\widehat{g_\lambda}(u)=(\Delta_{(-\delta_\lambda,
~\delta_\lambda)^m }~g_\lambda)^{^\wedge}(u)
$$
for almos all $u\in\mathbb R^m$. (We used the standard facts on
the convolution of functions in $L^1$ with functions in $L^p$, see
e.g. [11, Ch. I, \S~2].)

  Thus, according to Lemma 1, we see that if $\lambda>0$
is sufficiently large, then for almost all $u\in \lambda W$ we
have
$$
c_m\delta_\lambda^m\leq (\Delta_{(-\delta_\lambda,
~\delta_\lambda)^m }~g_\lambda)^{^\wedge}(u).
\eqno(9)
$$

   Since $\|\Delta_{(-\varepsilon, \varepsilon)^m}\|_{A(\mathbb R^m)}=
\Delta_{(-\varepsilon, \varepsilon)^m}(0)=1$, it follows that for
every function $f\in A_p(\mathbb R^m)$ and an arbitrary
$\varepsilon>0$ we have
$$
\|\Delta_{(-\epsilon, \epsilon)^m} f\|_{A_p(\mathbb R^m)}\leq
\|f\|_{A_p(\mathbb R^m)}.
$$
So, raising inequality (9) to the power $p$ and integrating over
$u\in\lambda W$, we see that
$$
(c_m^p\delta_\lambda^{mp}
|\lambda W|)^{1/p}\leq
\bigg(\int_{\lambda W}
|(\Delta_{(-\delta_\lambda,
~ \delta_\lambda)^m}g_\lambda)^{^\wedge}(u)|^p du \bigg)^{1/p}
$$
$$
\leq\|\Delta_{(-\delta_\lambda,
~\delta_\lambda)^m~}g_\lambda\|_{A_p(\mathbb R^m)}\leq
\|g_\lambda\|_{A_p(\mathbb R^m)}=
\|F_\lambda\|_{A_p(\mathbb R^m)}
$$
for all sufficiently large $\lambda$.

  Hence, taking (8) into account we obtain
$$
c_m \delta_\lambda^m\lambda^{m/p}|W|^{1/p}\leq
2\|e^{i\lambda\varphi}\|_{A_p(\mathbb R^m, V)}.
$$

   It remains only to note that condition (6) implies
$\delta_\lambda\geq c \chi^{-1}(1/\lambda)$. The theorem is
proved.

\quad

\emph{Remark.} Let $D$ be a bounded domain in $\mathbb R^n, ~n\geq
2$. Consider its characteristic function $1_D$, i.e., the function
that takes value $1_D(t)=1$ for $t\in D$ and value $1_D(t)=0$ for
$t\notin D$. One can show that if the boundary $\partial D$ of a
domain $D$ is $C^2$ -smooth, then $1_D\in A_p(\mathbb R^n)$ for
$p>\frac{2n}{n+1}$ and $1_D\notin A_p(\mathbb R^n)$ for
$p\leq\frac{2n}{n+1}$. At the same time in the more general case
of domains with $C^1$ -smooth boundary the Fourier transform of
the characteristic function may behave in an essentially different
way. The author constructed a domain $D\subset\mathbb R^2$ with
$C^1$ -smooth boundary such that $1_D\in A_p(\mathbb R^2)$ for all
$p>1$. Theorem $1''$ plays the key role in the study of the
question for which domains $D\subset\mathbb R^n$ with $C^1$
-smooth boundary we have inclusion $1_D\in A_p(\mathbb R^n)$. Our
results will be presented in another publication.

\quad

\begin{center}
\textbf{\S~2. Slow growth of $\|e^{i\lambda\varphi}\|_{A_p(\mathbb
T^m)}$}
\end{center}

  In the work [3] for each given class $C^{1,\omega}(\mathbb T)$
(under certain simple assumption imposed on $\omega$) we
constructed a nontrivial real function $\varphi\in
C^{1,\omega}(\mathbb T)$ with slow growth of the norms
$\|e^{i\lambda\varphi}\|_{A_p(\mathbb T)}$. Namely, put (as above
$\chi^{-1}$ is the function inverse to
$\chi(\delta)=\delta\omega(\delta)$)
$$
\Theta_1(y)=\frac{y}{\log y}\chi^{-1}\bigg(\frac{(\log y)^2}{y}\bigg)
$$
and for $1<p<2$ put
$$
\Theta_p(y)=\bigg(\int_1^y
\bigg(\chi^{-1}\bigg(\frac{1}{\tau}\bigg)\bigg)^p d\tau\bigg)^{1/p}.
$$
Let $\omega$ satisfies condition $\omega(2\delta)<2\omega(\delta)$
for all sufficiently small $\delta>0$. Then there exists a real
nowhere linear, i.e. not linear on any interval, function
$\varphi\in C^{1, \omega}(\mathbb T)$ such that for all $p, ~1\leq
p<2,$ we have $\|e^{i\lambda\varphi}\|_{A_p(\mathbb
T)}=O(\Theta_p(|\lambda|)), ~|\lambda|\rightarrow\infty,
~\lambda\in\mathbb R$.

  This result on the slow growth easily transfers to the case of
a torus of an arbitrary dimension. Let us say that the gradient of
a function $\varphi$ is nowhere degenerate if it is nondegenerate
on any open set, that is if for every open set $V\subseteq \mathbb
R^m$ we have $|\nabla \varphi (V)|>0$.

\quad

\textbf{Theorem 2.} \emph{Let $\omega(2\delta)<2\omega(\delta)$
for all sufficiently small $\delta>0$. There exists a real
function $\varphi\in C^{1, \omega}(\mathbb T^m)$ with nowhere
degenerate gradient such that for all $p, ~1\leq p<2,$ we have}
$$
\|e^{i\lambda\varphi}\|_{A_p(\mathbb
T^m)}=O((\Theta_p(|\lambda|))^m),
\qquad |\lambda|\rightarrow\infty, \quad \lambda\in\mathbb R.
$$

\quad

  The multidimensional version follows from the one-dimensional.
Indeed, if $\varphi_0\in C^{1,\omega}(\mathbb T)$ is a real
nowhere linear function, then, taking
$$
\varphi(t)=\varphi_0(t_1)+\varphi_0(t_2)+\ldots +\varphi_0(t_m),
\qquad t=(t_1, t_2, \ldots, t_m)\in\mathbb T^m,
$$
we have a function $\varphi\in C^{1,\omega}(\mathbb T^m)$ with
nowhere degenerate gradient and it remains only to note that
$\|e^{i\lambda\varphi}\|_{A_p(\mathbb
T^m)}=(\|e^{i\lambda\varphi_0}\|_{A_p(\mathbb T)})^m$.

   In the same manner (or directly from Theorem 2, taking
Corollary 1 into account) we obtain the corollaries below, which
are the multidimensional versions of Corollaries 2, 3 from our
paper [3].

\quad

\textbf{Corollary 3.} \emph{Let $0<\alpha<1$. There exists a real
function $\varphi\in C^{1, \alpha}(\mathbb T^m)$ with nowhere
degenerate gradient such that}
$$
\|e^{i\lambda\varphi}\|_{A(\mathbb T^m)}=
O(|\lambda|^{\frac{m\alpha}{1+\alpha}}
(\log |\lambda|)^{\frac{m(1-\alpha)}{1+\alpha}});
\leqno(\textrm{i})
$$
$$
\|e^{i\lambda\varphi}\|_{A_p(\mathbb T^m)}\simeq
|\lambda|^{m\big(\frac{1}{p}-\frac{1}{1+\alpha}\big)}
\qquad\qquad\mathit{for}\quad 1<p<1+\alpha,
\leqno(\textrm{ii})
$$
$$
\|e^{i\lambda\varphi}\|_{A_p(\mathbb T^m)}\simeq 1
\qquad\qquad\qquad\qquad\quad\mathit{for}\quad 1+\alpha<p<2,
$$
$$
\|e^{i\lambda\varphi}\|_{A_p(\mathbb T^m)}=O((\log |\lambda|)^{m/p})
\qquad\mathit{for}\quad p=1+\alpha.
$$

\quad

   In particular, we see that for $1<p<1+\alpha$ the estimate
in Corollary 1 of the present paper is sharp.

\quad

\textbf{Corollary 4.} \emph{Let $\gamma(\lambda)\geq 0$ and
$\gamma(\lambda)\rightarrow+\infty$ as $\lambda\rightarrow
+\infty$. There exists a real function $\varphi\in C^1(\mathbb
T^m)$ with nowhere degenerate gradient such that}
$$
\|e^{i\lambda\varphi}\|_{A(\mathbb T^m)}=
O(\gamma(|\lambda|)(\log |\lambda|)^m).
$$

\quad

\begin{center}
\textbf{\S~3. Upper estimates}
\end{center}

  Let $1\leq p<2$. As we noted in Introduction, if
$\varphi$ is a real function on the circle $\mathbb T$ satisfying
Lipschitz condition of order 1, then
$\|e^{i\lambda\varphi}\|_{A_p(\mathbb T)}=O(|\lambda|^{1/p-1/2})$
and due to Leibenson--Kahane--Alp\'ar estimate (2) we have
$\|e^{i\lambda\varphi}\|_{A_p(\mathbb
T)}\simeq|\lambda|^{1/p-1/2}$ for every nonlinear real function
$\varphi\in C^2(\mathbb T)$ (actually this holds even under
assumption that $\varphi\in C^{1,1}(\mathbb T)$, see Corollary 2
for $m=1$). Here we shall obtain similar results in the case of
torus of an arbitrary dimension.

   Recall the known smoothness condition that guarantees
that a function belongs to the classes $A_p(\mathbb T^m), ~1\leq
p< 2$, namely: if $f\in C^{\nu, \alpha}(\mathbb T^m)$ and
$\nu+\alpha>m(1/p-1/2)$, then $f\in A_p(\mathbb T^m)$. It is known
that this condition is sharp in the sense that for
$\nu+\alpha=m(1/p-1/2)$ the inclusion $C^{\nu, \alpha}(\mathbb
T^m)\subseteq A_p(\mathbb T^m)$ fails. \footnote{In the
one-dimensional case for $p=1$ the corresponding results are due
to Bernstein (see [12, Ch. VI, Theorem 3.1]), the generalization
to the case $1<p<2$ was obtained by Sz\'asz (see [12, Ch. VI,
Theorem 3.10]), in the multidimensional case the sufficiency of
the indicated condition was obtained by Sz\'asz and
Minakshisundaram [13] (see also remark of Bochner [14]). The
sharpness of the indicated smoothness condition in the
multidimensional case was established by Weinger [15].}

  We shall show that the following theorem holds.

\quad

   \textbf{Theorem 3.} \emph{Let $1\leq p<2$ and let $0\leq\alpha\leq 1$.
Let $\varphi\in C^{\nu, \alpha}(\mathbb T^m)$ be a real function.
Suppose that $\nu+\alpha\geq 1, ~\nu+\alpha>m(1/p-1/2)$. Then}
$$
\|e^{i\lambda\varphi}\|_{A_p(\mathbb T^m)}= O(|\lambda|^{m(1/p-1/2)}),
\qquad |\lambda|\rightarrow\infty,  \quad \lambda\in\mathbb R.
$$

\quad

  Under assumption that $\nu+\alpha\geq 2$, we have
$C^{\nu, \alpha}\subseteq C^{1,1}$, so, using Corollary 2 and
Theorem 3, we see that the following theorem holds.

\quad

  \textbf{Theorem 4.} \emph{Let $1\leq p<2$ and let $0\leq\alpha\leq 1$.
Let $\varphi\in C^{\nu, \alpha}(\mathbb T^m)$ be a real function
with nondegenerate gradient. Suppose that $\nu+\alpha\geq 2,
~\nu+\alpha>m(1/p-1/2)$. Then
$$
\|e^{i\lambda\varphi}\|_{A_p(\mathbb T^m)}\simeq|\lambda|^{m(1/p-1/2)},
\qquad |\lambda|\rightarrow\infty, \quad \lambda\in\mathbb R.
$$
In particular, if $\nu+\alpha\geq 2, ~\nu+\alpha>m/2$, then
$\|e^{i\lambda\varphi}\|_{A(\mathbb T^m)}\simeq |\lambda|^{m/2}$.}

\quad

  \emph{Proof of Theorem} 3. In the usual way we define
the norm on the space $L^2(\mathbb T^m)$ and the norm on the space
$C(\mathbb T^m)$ of continuous functions on $\mathbb T^m$:
$$
\|f\|_{L^2(\mathbb T^m)}=
\bigg(\frac{1}{(2\pi)^m}\int_{\mathbb T^m}|f(t)|^2 dt\bigg)^{1/2},
\qquad \|g\|_{C(\mathbb T^m)}=\sup_{t\in\mathbb T^m}|g(t)|.
$$

   We define the norm on $C^{\nu, \alpha}(\mathbb T^m)$ by
$$
\|f\|_{C^{\nu, \alpha}(\mathbb T^m)}=\max_{0\leq \gamma_1+\gamma_2+\ldots
+\gamma_m\leq \nu}
\|D_{\gamma_1, \gamma_2, \ldots, \gamma_m}f\|_{C(\mathbb T^m)}+
$$
$$
+\max_{\gamma_1+\gamma_2+\ldots +\gamma_m=\nu}
\sup_{\delta>0}\frac{1}{\delta^\alpha}
\omega(D_{\gamma_1, \gamma_2, \ldots, \gamma_m}f, \delta),
$$
where
$$
D_{\gamma_1, \gamma_2, \ldots, \gamma_m}f(t)=
\frac{\partial^{\gamma_1+\gamma_2+\ldots +\gamma_m} f}{\partial
t_1^{\gamma_1}\partial t_2^{\gamma_2}\ldots \partial
t_m^{\gamma_m}}, \qquad t=(t_1, t_2, \ldots, t_m).
$$

   It is easy to see that for $\nu=0, 1, 2, \ldots$
and $0\leq \alpha\leq 1$ we have
$$
\|fg\|_{C^{\nu, \alpha}(\mathbb T^m)}\leq c \|f\|_{C^{\nu,
\alpha}(\mathbb T^m)}\|g\|_{C^{\nu, \alpha}(\mathbb T^m)},
$$
where $c=c(m, \nu, \alpha)>0$ is independent of $f$ and $g$.

  Using this relation it is easy to verify that if
$\nu+\alpha\geq 1$, then for each real function $\varphi\in
C^{\nu, \alpha}(\mathbb T^m)$ we have the estimate
$$
\|e^{i\lambda\varphi}\|_{C^{\nu, \alpha}(\mathbb T^m)} \leq c
|\lambda|^{\nu+\alpha}, \qquad \lambda\in\mathbb R, \quad |\lambda|\geq 1,
$$
where $c=c(\nu, \alpha, \varphi)$ is independent of $\lambda$. For
$\nu=1, 2, \ldots$ this can be verified by induction over $\nu$
with fixed $\alpha$. For $\nu=0$ the condition $\nu+\alpha\geq 1$
yields $\alpha=1$ and the required estimate also holds.

  Thus it is clear that the statement of the theorem is an
immediate consequence of the following simple lemma.

\quad

\textbf{Lemma 2.} \emph{Let $f\in C^{\nu, \alpha}(\mathbb T^m)$
and let $\nu+\alpha>m(1/p-1/2), ~1\leq p<2$. Then
$$
\|f\|_{A_p(\mathbb T^m)}\leq c \|f\|_{C^{\nu, \alpha}(\mathbb
T^m)}^\tau\|f\|_{L^2(\mathbb T^m)}^{1-\tau},
$$
where
$$
\tau=\frac{m(1/p-1/2)}{\nu+\alpha}.
$$}

\quad

\emph{Proof.} Let $S^{m-1}$ denote the unit sphere in $\mathbb
R^m$ centered at $0$. Let $\delta>0$ and let $\xi\in S^{m-1}$. For
each $j=1, 2, \ldots , m$ consider a function
$$
\frac{\partial^\nu f}{\partial t_j^\nu}(t+\delta\xi)-\frac{\partial^\nu
f}{\partial t_j^\nu}(t-\delta\xi), \qquad t=(t_1, t_2, \ldots , t_m)\in
\mathbb T^m.
$$
Writing the Parseval identity for this function we obtain
$$
\frac{1}{(2\pi)^m}\int_{\mathbb T^m}\bigg| \frac{\partial^\nu f}{\partial
t_j^\nu}(t+\delta\xi)-\frac{\partial^\nu f}{\partial
t_j^\nu}(t-\delta\xi) \bigg|^2 dt =
\sum_{k=(k_1, k_2, \ldots, k_m)\in\mathbb Z^m}
k_j^{2\nu}|\widehat{f}(k)|^2 4\sin^2 (\delta k, \xi).
$$
Thus,
$$
\sum_{k\in\mathbb Z^m} k_j^{2\nu}|\widehat{f}(k)|^2 \sin^2 (\delta k,
\xi)\leq \|f\|_{C^{\nu, \alpha}}^2 \delta^{2\alpha}.
$$
Summing over $j=1, 2, \ldots , m$, we have
$$
\sum_{k\in\mathbb Z^m} |k|^{2\nu}|\widehat{f}(k)|^2 \sin^2 (\delta k,
\xi)\leq c\|f\|_{C^{\nu, \alpha}}^2 \delta^{2\alpha},
$$
where $c=c(\nu, m)>0$.

   Let $|k|\leq 1/\delta$, then $|(\delta k, \xi)|\leq 1$ whence
$|\sin (\delta k, \xi)|\geq |(\delta k, \xi)|/2 $ and we see that
$$
\sum_{|k|\leq 1/\delta} |k|^{2\nu}|\widehat{f}(k)|^2 (\delta k,
\xi)^2\leq  c\|f\|_{C^{\nu, \alpha}}^2 \delta^{2\alpha}.
$$
Integrating this inequality over $\xi\in S^{m-1}$ (for $m=1$ we
just put $\xi=1$) and taking into account the fact that for an
arbitrary vector $v\in\mathbb R^m$
$$
\int_{S^{m-1}} (v, \xi)^2 d\xi= c_m |v|^2,
$$
we obtain
$$
\sum_{|k|\leq 1/\delta} |k|^{2\nu}|\widehat{f}(k)|^2 \delta^2 |k|^2\leq
c\|f\|_{C^{\nu, \alpha}}^2 \delta^{2\alpha}.
$$

   Let $\delta=2^{-n}, ~n=1, 2, \ldots$. We see that
$$
\sum_{2^{n-1}\leq |k|<2^n}|\widehat{f}(k)|^2\leq c \|f\|_{C^{\nu,
\alpha}}^2 2^{-2n(\nu+\alpha)}, \qquad n=1, 2, \ldots ,
\eqno(10)
$$
with constant $c>0$ independent of $n$ and $f$.

   Using the H\"{o}lder inequality with $p^*=2/p, ~1/p^*+1/q^*=1,$
we obtain from relation (10) that
$$
\sum_{2^{n-1}\leq |k|<2^n}|\widehat{f}(k)|^p\leq \bigg (\sum_{2^{n-1}\leq
|k|<2^n}|\widehat{f}(k)|^{pp^*}\bigg )^{1/p^*} \bigg ( \sum_{2^{n-1}\leq
|k|<2^n} 1 \bigg )^{1/q^*}
$$
$$
\leq c \|f\|_{C^{\nu, \alpha}}^p 2^{-pn(\nu+\alpha)}(2^{nm})^{1-p/2}=
c \|f\|_{C^{\nu, \alpha}}^p 2^{-np(\nu+\alpha)(1-\tau)}.
$$
Hence it is clear that for an arbitrary $B\geq 1$ we have
$$
\sum_{|k|\geq B}|\widehat{f}(k)|^p\leq c \|f\|_{C^{\nu, \alpha}}^p
B^{-p(\nu+\alpha)(1-\tau)}.
\eqno(11)
$$

  At the same time it is obvious that for $B\geq 1$
$$
\sum_{|k|<B}|\widehat{f}(k)|^p\leq \bigg
(\sum_{|k|<B}|\widehat{f}(k)|^{pp^*}\bigg )^{1/p^*}
\bigg(\sum_{|k|<B} 1\bigg )^{{1/q^*}}
$$
$$
\leq\|f\|_{L^2}^p c B^{m(1-p/2)}=c \|f\|_{L^2}^p
B^{p(\nu+\alpha)\tau}.
\eqno(12)
$$

  It remains to add relations (11), (12) and put
$$
B=\bigg (\frac{\|f\|_{C^{\nu, \alpha}}}
{\|f\|_{L^2}}\bigg )^{\frac{1}{\nu+\alpha}}.
$$
The lemma is proved. Theorem 3 is proved. Theorem 4 follows.

\quad

\quad

\quad

\quad

\begin{center}
\textbf{References}
\end{center}
\flushleft
\begin{enumerate}

\item A. Beurling, H. Helson, ``Fourier--Stieltjes transforms
    with bounded powers'', \emph{Math. Scand.}, \textbf{1}
    (1953), 120--126.

\item J.-P. Kahane, \emph{S\'{e}ries de Fourier absolument
    convergentes}, Springer-Verlag, Berlin - Heidelberg - New
    York, 1970.

\item  V. V. Lebedev, ``Quantitative estimates in Beurling --
    Helson type theorems'', \emph{Sbornik: Mathematics},
    \textbf{201}:12, (2010), 1811--1836.

\item  Z. L. Leibenson, ``On the ring of functions with
    absolutely convergent Fourier series'', \emph{Usp. Mat.
    Nauk.}, \textbf{9}:3(61) (1954), 157--162 (in Russian).

\item  J.-P. Kahane, ``Sur certaines classes de s\'eries de
    Fourier absolument convergentes'', \emph{J. de
    Math\'ematiques Pures et Appliqu\'ees}, \textbf{35}:3
    (1956), 249--259.

\item  L. Alp\'ar, ``Sur une classe partiquli\`ere de s\'eries
    de Fourier \`a certaines puissances absolument
    convergentes'', \emph{Studia Sci. Math. Hungarica},
    \textbf{3} (1968), 279--286.

\item V. V. Lebedev, ``Diffeomorphisms of the circle and the
    Beurling -- Helson theorem'', \emph{Functional analysis
    and its applications}, \textbf{36}:1 (2002), 25--29.

\item  G. W. Hedstrom, ``Norms of powers of absolutely
    convergent Fourier series in several variables'',
    \emph{Michigan Math, J.}, \textbf{14}:4 (1967), 493-495.

\item  M. V. Korobkov, ``Properties of $C^1$ -smooth functions
    whose gradient range has topological dimension 1'',
    \emph{Doklady Mathematics (Doklady Akademii Nauk)},
    \textbf{81}:1 (2010), 11--13.

\item  M. Plancherel, G. P\'olya, ``Fonctions enti\`eres et
    int\'egrales de Fourier multiples. II'', \emph{Comment.
    Math. Helv.}, \textbf{10}:2 (1937), 110-163.

\item  E. M. Stein and G. Weiss, \emph{Introduction to Fourier
    analysis on Euclidean spaces}, Princeton Univ. Press,
    Princeton, NJ, 1971.

\item  A. Zigmund, \emph{Trigonometric series}, vol. I,
    Cambrige Univ. Press, New York, 1959.

\item  S. Minakshisundaram, O. Sz\'asz, ``On absolute
    convergence of multiple Fourier series'', \emph{Trans.
    Amer. Math. Soc.}, \textbf{61}:1 (1947), 36--53.

\item  S. Bochner, ``Reviews of `On absolute convergence of
    multiple Fourier series' by Sz\'asz and
    Minakshisundaram'', \emph{ Math. Rev.}, \textbf{8} (1947),
    376.

\item  S. Weinger, ``Special trigonometric series in $k$
    dimensions'', \emph{Mem. Amer. Math. Soc.}, \textbf{59}
    (1965).

\end{enumerate}

\quad

\qquad\textsc{V. V. Lebedev}\\
\qquad Dept. of Mathematical Analysis\\
\qquad Moscow State Institute of\\
\qquad Electronics and Mathematics (Technical University)\\
\qquad E-mail address: \emph {lebedevhome@gmail.com}

\end{document}